\documentclass{amsart}
\usepackage[cp1251]{inputenc}
\usepackage[english]{babel}
\usepackage{amsmath}
\usepackage{amssymb}
\usepackage{amsfonts}
\newtheorem{lem}{Lemma}

\newtheorem{cor}{Corollary}
\begin{document}

\thispagestyle{empty}

 \title[Modeling of conservative fields]{Modeling of conservative fields in piecewise infinite plate with a
thin inclusion
 }%

\author{O.Yaremko;E. Mogileva}%Указываем авторов

\address{Oleg Emanuilovich Yaremko,
\newline\hphantom{iii}Penza State University,% Место работы
\newline\hphantom{iii}str. Lermontov, 37, % Адрес (улица, дом, строение и т.п.)
\newline\hphantom{iii} 440038, Penza, Russia}%  Адрес (почтовый индекс, город, страна)
\email{yaremki@mail.ru}% Ваш электронный адрес для переписки

\maketitle {\small

\begin{quote}
\noindent{\sc Abstract. }
In the article the formulas for the modeling of
conservative fields in piecewise infinite plate with a thin inclusion found.
The accuracy of the found formulas is of order equal to the thickness of the
outer layer. The problem for higher accuracy solved by means of asymptotic
formulas.

\medskip


 \end{quote}

\section{Introduction}

Problems about structure of the conservative field of a two-layer flat plate
leads to a separate system of the Laplace's equation with a boundary
condition to Dirichlet and internal coupling conditions [1], [2]. In the
work the method of transformation operators is applied. Yaremko O. E.
developed this method [1], [2], [3], [4], [5]. As a result the conservative
field for the two-layer piecewise plate can be interpreted as a deformation
of the conservative field of a homogeneous piecewise plate. Thus the
deforming transformation operator is written out. The deforming
transformation operator is recursive, and convenient for practical
realization on the computer. In this case, use the analytic representation
of a field by means of the deforming transformation operator is difficult.

We consider the case of small thickness l cross-border layer of the plate.
Let's consider that the physical properties of the layers very different,
that is the ratio of the thermal conductivity of the layers close to zero.

The method of the deforming transformation operator gives algorithm
inconvenient for realization. Formulas which have no these defect are found
based on the Euler--Maclaurin formula.

Numerical methods of research lead to ill-conditioned [8] systems of the
linear algebraic equations.

\subsection{ The half-plane with internal coupling condition}

Modelling of potential fields for a semi-limited infinite plate with thin
inclusion leads to the Dirichlet problem for two-layer half-plane with
coupling conditions on the direct. Let's consider a boundary problem about
the solving of separate system of the Laplace's equation
\begin{equation}
\label{eq1}
a_1^2 \frac{\partial ^2u_1 }{\partial x^2}+\frac{\partial ^2u_1 }{\partial
y^2}=0,(x,y);0<x<l,-\infty <y<\infty
\end{equation}
\begin{equation}
\label{eq2}
a_2^2 \frac{\partial ^2u_2 }{\partial x^2}+\frac{\partial ^2u_1 }{\partial
y^2}=0,(x,y);l<x,-\infty <y<\infty
\end{equation}
by the boundary condition for direct $x=l:$
\begin{equation}
\label{eq3}
u_1 (x,y)=\hat {u}(x,y),x=0
\end{equation}
by the internal coupling conditions [1], [2], [3], [4], [5] for direct $x=l:$
\begin{equation}
\label{eq4}
u_1 (x,y)=u_2 (x,y),k\frac{\partial u_1 }{\partial x}(x,y)=\frac{\partial
u_2 }{\partial x}(x,y),k>0,
\end{equation}
Here function $\hat {u}=\hat {u}(x,y)$- is harmonic in the right half-plane,
$H=\{(x,y):0<x,y\in R\}$ and continuous in $\{(x,y):0\le x,y\in R\}$ and the
conditions are met
\[
\hat {u}=0(\sqrt {x^2+y^2} ),
\]
\[
\int\limits_{-\infty }^\infty {\frac{\hat {u}(0,y)}{1+y^2}dy<\infty .}
\]
The problem (\ref{eq1})-(\ref{eq4}) is ill-posed [8] in the case of a thin inclusion, that
is, small value l , and in the case of a strong heterogeneity, that is,
large differences in the parameters $a_1 ,a_2 $. The solution of this
problem by known methods is accompanied by considerable computing
difficulties. New analytical methods to solve this problem it is necessary
to develop.

\subsection{Thin strip }

Function $\hat {u}=\hat {u}(x,y)$ is a harmonic function in the right
half-plane $H=\{(x,y):0<x,y\in R\}$ continuous in $\{(x,y):0\le x,y\in R\}$
and the conditions are met
\[
\hat {u}=0(\sqrt {x^2+y^2} ),
\]
\[
\int\limits_{-\infty }^\infty {\frac{\hat {u}(0,y)}{1+y^2}dy<\infty .}
\]
Let's consider the Dirichlet problem for Laplace's equation in the strip
$H_l =\{(x,y):0<x<l,y\in R\}$
\begin{equation}
\label{eq5}
\frac{\partial ^2u_1 }{\partial x^2}+\frac{\partial ^2u_1 }{\partial
y^2}=0,(x,y);0<x<l,-\infty <y<\infty
\end{equation}
with boundary conditions
\begin{equation}
\label{eq6}
u(0,y)=\hat {u}(0,y),u(l,y)=0.
\end{equation}
The problem (\ref{eq5})-(\ref{eq6}) is ill-posed [8] at small value l.

\subsection{The circle with the internal coupling conditions }

Let the function $\hat {u}=\hat {u}(x,y)$- is a harmonic function in the
unit circle, $B=\{(x,y):x^2+y^2<1\}$ and continuous in its closure $\bar
{B}$. The boundary problem about the solving a separate system of the
Laplace's equation consider in the circle B
\[
\Delta u_1 =0,(x,y)\in K_R ,
\]
\begin{equation}
\label{eq7}
\Delta u_2 =0,(x,y)\in B_R
\end{equation}
Here $K_R $ - annulus with the smaller radius R and larger radius 1, $B_R $
- the circle of radius R. Let the boundary condition on the circle satisfied
S:
\begin{equation}
\label{eq8}
u_1 (x,y)=\hat {u}(x,y),(x,y)\in S
\end{equation}
internal coupling condition on a circle of radius R:
\[
u_1 (x,y)=u_2 (x,y),
\]
\begin{equation}
\label{eq9}
kL_0 u_1 (x,y)=L_0 u_2 (x,y),k>0,(x,y)\in S_R ,
\end{equation}
\[
L_0 =x\frac{\partial }{\partial x}+y\frac{\partial }{\partial y}.
\]
Let's notice that the problem (\ref{eq7})-(\ref{eq9}) is ill-posed [8] at values $R\approx
1.$

\subsection{ Thin annulus}

Let the function $\hat {u}=\hat {u}(x,y)$ - is a harmonic function in the
unit circle, $B=\{(x,y):x^2+y^2<1\}$ and continuous in its closure $\bar
{B}$. In annulus $K_R $- with the smaller radius R and larger radius 1,
consider the Dirichlet problem about the solution of the Laplace's equation
\begin{equation}
\label{eq10}
\Delta u_1 =0,(x,y)\in K_R ,
\end{equation}
by the boundary conditions on the circumferences S, $S_R $:
\begin{equation}
\label{eq11}
u=\hat {u},(x,y)\in S;u=0,(x,y)\in S_R .
\end{equation}
The problem (\ref{eq10}) is ill-posed [8] by the values $R\approx 1.$

\section{Transformation operators}

In the works of the author [1], [2], [3], [4] and [5] the concept of
operator transformation J, $J:\hat {u}\to u,$ allowing a known solution
$\hat {u}$ of the model problem to solve u any of the four above-listed
problems determined.

\subsection{The half-plane with the internal coupling conditions }

Transformation operator J has the form
\[
J:\hat {u}\to u,u=\chi (H_1 )u_1 +\chi (H_2 )u_2 ,
\]
where:
\[
u_1 (x,y)=\sum\limits_{j=0}^\infty {\left( {\frac{1-k}{1+k}} \right)^j\left(
{\hat {u}(x+2lj,y)-\frac{1-k}{1+k}\hat {u}(2l-x+2lj,y)} \right)} ,
\quad
0<x<l,
\]
\begin{equation}
\label{eq12}
u_2 (x,y)=\frac{2k}{k+1}\sum\limits_{j=0}^\infty {\left( {\frac{1-k}{1+k}}
\right)^j\hat {u}\left( {\frac{a_1 }{a_2 }(x-l)+l+2lj,y} \right)} ,
\quad
l<x,
\end{equation}
\[
k=\frac{\lambda _1 }{\lambda _2 }\frac{a_2 }{a_1 },
\]
$\chi -$ is the characteristic function of the set [6].

\subsection{Thin strip}

Transformation operator $J:\hat {u}\to u,$ can be determined by the formula
\begin{equation}
\label{eq13}
u(x,y)=\sum\limits_{j=0}^\infty {\left( {\hat {u}(x+2lj,y)-\hat
{u}(2l-x+2lj,y)} \right)} ;
\quad
0<x<l.
\end{equation}
\subsection{ Circle with the internal coupling conditions}

Transformation operator is expressed by the formula:
\[
J:\hat {u}\to u,u=\chi (K_R )u_1 +\chi (B_R )u_2 ,
\]
where:
\begin{equation}
\label{eq14}
u_1 (x,y)=\sum\limits_{j=0}^\infty {\left( {\frac{1-k}{1+k}} \right)^j\left(
{\hat {u}(xR^{2lj},yR^{2lj})-\frac{1-k}{1+k}\hat {u}\left(
{\frac{R^{2j+2}}{x},\frac{R^{2j+2}}{y}} \right)} \right)} ,
\end{equation}
\[u_2 (x,y)=\frac{2k}{k+1}\sum\limits_{j=0}^\infty {\left( {\frac{1-k}{1+k}}
\right)^j\hat {u}(xR^{2j},yR^{2j})} ,
\]

$\chi \quad -$ is the characteristic function of the set [6].

\subsection{ Thin annulus}

Transformation operator $J:\hat {u}\to u,$ can be determined by the formula
\begin{equation}
\label{eq15}
u(x,y)=\sum\limits_{j=0}^\infty {\left( {u\left( {xR^{2j},yR^{2j}}
\right)-u\left( {\frac{R^{2j+2}}{x},\frac{R^{2j+2}}{y}} \right)} \right)} ,
\end{equation}
Formulas (\ref{eq12}) and (\ref{eq14}) are useful only if the values $k\approx 1.$ In the
case of the two-layer body, which components are sharply different
properties, that is at values k close to zero, and for large values k,
transformation operators cannot be used, in view of the slow convergence of
the series at $\frac{1-k}{1+k}\approx 1.$

Formulas(\ref{eq13}) and (\ref{eq15}) cannot be used in view of the slow convergence of the
series at $l\approx 1$ and $R\approx 1$, respectively.

The purpose of our research is to consists in creation of a design of
transformation operators comfortable in these cases.

\section{ Auxiliary result}

\begin{lem} If the function y = f(x) is defined on the line segment
[0,1] and has limited variation [10] on the line segment [0,1], then for
each value R the following assessment for the difference of integral and its
integral sum is carried
\[
\left| {\int\limits_0^1 {x^{-1}f(x)dx} -\ln
\frac{1}{R^2}\sum\limits_{j=0}^\infty f(R^{2j})} \right|\le \ln
\frac{1}{R^2}\mathop V\limits_0^1 (f)
\]
\end{lem}
\textbf{Proof}.Let's spread integral into the sum (n +1) of the summand:

$\int\limits_0^1 {x^{-1}f(x)dx} =\sum\limits_{j=0}^{n-1} {\int\limits_{R^2}^1
{x^{-1}f(R^{2j}x)dx+\int\limits_0^1 {x^{-1}f(R^{2n}x)dx} } } $Variable
substitution in the last integral results in equality:
\[
\int\limits_0^1 {x^{-1}f(x)dx} =\sum\limits_{j=0}^{n-1} {\int\limits_{R^2}^1
{f(R^{2j})\frac{dx}{x}+\int\limits_0^{R^{2n}} {x^{-1}f(x)dx} } }
\]
Considering that $$\left| {f(R^{2j}x)-f(R^{2j})} \right|\le \mathop
V\limits_0^1 (f),$$ We obtain the required assessment:
\[
\left| {\int\limits_0^1 {x^{-1}f(x)dx-\sum\limits_{j=0}^\infty
{\int\limits_{R^2}^1 {f(R^{2j})\frac{dx}{x}} } } } \right|\le
\]
\[
\le \sum\limits_{j=0}^\infty {\int\limits_{R^2}^1 {\left|
{f(R^{2j}x)-f(R^{2j})} \right|\frac{dx}{x}\le \mathop V\limits_0^1 (f)\ln
\frac{1}{R^2}} }
\]
$\mathop V\limits_0^1 (f)\quad -$ variation of function [10] f on the line
segment [0,1]. The flat analog of the Lemma is well known [8].
\begin{lem} If the function $y = f(x)$ is defined on the ray $[0;\infty
)$ and has limited variation [10] on the ray $[0; \infty )$; then for each
value l the following assessment for the difference of integral and its
integral sum is carried
\[
\left| {\int\limits_0^\infty {f(x)dx-2l\sum\limits_{j=0}^\infty {f(2lj)} } }
\right|\le 2l\mathop V\limits_0^\infty (f)
\]
The Bernoulli numbers [11] are defined using the generating function
\[
\frac{z}{e^z-1}=\sum\limits_{j=0}^\infty {\frac{B_j }{j!}z^n.}
\]
The Euler--Maclaurin formula applied to the function $f(2lx)$ on the ray
$[0;\infty )$ results in equality

$$\sum\limits_{j=0}^\infty {f(2lj)\cong \frac{1}{2l}\int\limits_0^\infty
{f(x)dx} +\frac{f(0)}{2}} -\sum\limits_{k=1}^\infty {\frac{(2l)^{2k-1}B_{2k}
}{k!}f^{2k-1}(0).}
$$
$$\sum\limits_{j=0}^\infty {f(R^{2j})\cong \frac{1}{\ln
\frac{1}{R^2}}\int\limits_0^\infty {f(x)dx} +\frac{f(\ref{eq1})}{2}}
-\sum\limits_{k=1}^\infty {\frac{(2l)^{2k-1}B_{2k} }{k!}L_0^{2k-1} f(\ref{eq1})}
$$
\end{lem}
\begin{cor} When $k\in (0,1)$ completed asymptotically
\[
\sum\limits_{j=0}^\infty {\left( {\frac{1-k}{1+k}} \right)^jf(x+2lj)\cong }
\]
$$\cong \frac{1}{2l}\int\limits_0^\infty {e^{h\varepsilon }f(x+\varepsilon
)d\varepsilon +\frac{f(x)}{2}} -\sum\limits_{k=1}^\infty
{\frac{(2l)^{2k-1}B_{2k} }{k!}L_h^{2k-1} f(x).} $$

When $k\in (1,\infty
)$completed asymptotically
$$\sum\limits_{j=0}^\infty {\left( {\frac{1-k}{1+k}} \right)^jf(x+2lj)\cong }
\frac{f(x)}{2}-\sum\limits_{k=1}^\infty {\frac{(2l)^{2k-1}B_{2k}
}{k!}(2^{2k}-1)L_h^{2k-1} f(x).} $$
The operator $L_{2h} $ has the form $L_{2h}
=2h+\frac{d}{dx}.$
\end{cor}
\begin{cor} When $k\in (0,1)$ completed asymptotically
\[
\sum\limits_{j=0}^\infty {\left( {\frac{1-k}{1+k}} \right)^jf(rR^{2j})\cong
}
\]
$$\cong \frac{1}{\ln \frac{1}{R^2}}\int\limits_0^1 {\varepsilon
^{h-1}f(r\varepsilon )d\varepsilon }
+\frac{f(r)}{2}+\sum\limits_{k=1}^\infty {\frac{\left( {\ln \frac{1}{R^2}}
\right)^{2k-1}B_{2k} }{k!}L_{2h}^{2k-1} f(r).} $$
When $k\in (1,\infty )$
completed asymptotically

$$\sum\limits_{j=0}^\infty {\left( {\frac{1-k}{1+k}} \right)^jf(rR^{2j})\cong
} \frac{f(r)}{2}+\sum\limits_{k=1}^\infty {\frac{(\ln
\frac{1}{R^2})^{2k-1}B_{2k} }{k!}(2^{2k}-1)L_{2h}^{2k-1} f(r).} $$
The operator
$L_{2h} $ has the form $L_{2h} =2h+\frac{d}{dx}.$
\end{cor}
\section{Main result}

\subsection{The half-plane with the internal coupling conditions}

Study the case of a thin shell, that is the case of a small thickness of the
final layer l. Physically this means that the coefficients of heat capacity
layers are very different.

\textbf{Theorem 1.} Let the condition is fulfilled
\[
\frac{1-k}{1+k}=e^{2hl},
\]
and the value k is small. For component $u_1 $ and $u_2 $ of the problem
solution (\ref{eq1})-(\ref{eq4}) approximate formulas are valid
\[
u_2 \approx \frac{1-e^{2hl}}{2l}\hat {u}_3 ,
\]
\[
u_1 \approx \frac{1}{2l}\hat {u}_3 \left( {x,y}
\right)-\frac{e^{2h}}{2l}\hat {u}_3 (2l-x,y),(x,y)\in H_2 ,
\]
Where $\hat {u}_3 $ - solution of the problem of Robin [12] with the
boundary condition
\[
2h\hat {u}_3 +\frac{\partial }{\partial n}\hat {u}{ }_3=\hat {u},(x,y),x=0.
\]
The assessment is fair
$$\left| {\frac{1-e^{2hl}}{2l}\int\limits_0^\infty {e^{\varepsilon h}\hat {u}}
(x+\varepsilon ,y)d\varepsilon -u_2 (x,y)} \right|\le (1-e^{2hl})\mathop
V\limits_0^\infty (e^{\varepsilon h}\hat {u}(x+\varepsilon
,y)).$$
\textbf{Proof}. Let's apply the lemma 2. As a result we obtain the
formula:
$$\left| {\frac{2k}{k+1}\frac{1}{2l}\int\limits_0^\infty {e^{\varepsilon
h}\hat {u}(x+\varepsilon ,y)d\varepsilon -u_2 (x,y)} } \right|\le
\frac{2k}{k+1}\mathop V\limits_0^\infty (e^{\varepsilon h}\hat
{u}(x+\varepsilon ,y)).$$
We will express k through h. Get the assessment,
which proved for $u_2 $. The formula links solutions the Dirichlet problem
and the Robin problem from [12]
\[
u_2 =-\int\limits_0^\infty {e^{\varepsilon h}\hat {u}(x+\varepsilon
,y)d\varepsilon .}
\]
This theorem allows the physical interpretation: component of the solution
$u_2 $ with the accuracy to a numerical multiplier it is approximately equal
to the solution of the third homogeneous boundary problem with boundary
condition
\[
2h\hat {u}_3 +\frac{\partial }{\partial n}\hat {u}{ }_3=\hat {u},(x,y),x=0.
\]
Similarly

\textbf{Theorem 2.} Let the condition is fulfilled
\[
\frac{1-k}{1+k}=e^{2hl},
\]
and the value $k>>1.$ For the problem solution (\ref{eq1})-(\ref{eq4}) approximate
formulas are valid
\[
u_1 \approx \frac{1-e^{2hl}}{1+e^{2hl}}\frac{1}{2l}\left( {\left( {\hat
{u}_3 (x,y)-e^{2hl}\hat {u}_3 (x+2l,y)} \right)+...} \right),
\]
\[
\left( {...+e^{2hl}(\hat {u}_3 (2l-x,y)-e^{2hl}\hat {u}_3 (4l-x,y))}
\right),
\]
\[
u_2 \approx \frac{1-e^{2hl}}{2l}\frac{1}{2l}\left( {\hat {u}_3
(x,y)-e^{2hl}\hat {u}_3 (x+2l,y)} \right).
\]
Theorem 2 is proved on the basis of the Lemma 2.

\subsection{ Thin strip}

We apply Lemma 2. For the solution of the Dirichlet problem in the strip
(\ref{eq5})-(\ref{eq6}) assessment is fair
\[
\left| {\int\limits_0^\infty {\frac{\hat {u}(x+\varepsilon ,y)-\hat
{u}(2l-x+\varepsilon ,y)}{2l}d\varepsilon -u(x,y)} } \right|\le \mathop
V\limits_0^\infty ,
\]
here
\[
\mathop V\limits_0^\infty =\mathop V\limits_0^\infty \left( {\hat
{u}(x+\varepsilon ,y)-\hat {u}(2l-x+\varepsilon ,y)} \right)
\]
- variation of function [10] taken from a variable $\varepsilon $ on the
interval [0;$\infty )$. Thus, the following theorem is obtained

\textbf{Theorem 3.} The solution of the Dirichlet problem in the strip
(\ref{eq5})-(\ref{eq6}) in the strip can be found on the approximate formula
\[
u(x,y)\approx \frac{\hat {u}_2 (x,y)-\hat {u}_2 (2l-x,y)}{2l},
\]
Where $\hat {u}{ }_2$ - the solution of the Neumann problem for the equation
Laplace's in a circle with the boundary condition
\[
\frac{\partial \hat {u}_2 }{\partial x}=\hat {u},x=0.
\]
\subsection{ The axial case}

The circle with with the conditions internal coupling. Study the case of a
thin shell, that is the case of a small thickness of the outer layer 1-R.
Let's consider also magnitude k of the small, that is consider the case when
the coefficients of heat capacity layers differ greatly. Assuming
\[
\frac{1-k}{1+k}=R^{2h}.
\]
Formula follow from Lemma 1

$$\left| {\frac{2k}{k+1}\frac{1}{\ln \frac{1}{R^2}}\int_0^1
{\varepsilon ^{h-1}\hat {u}(x\varepsilon ,y\varepsilon )d\varepsilon -u_2
(x,y)} } \right|\le$$
$$\le \frac{2k}{k+1} V_0^1  {\left(
{\varepsilon ^h\hat {u}(x\varepsilon ,y\varepsilon )} \right)},
$$
$$\left| {\frac{1-R^{2h}}{\ln \frac{1}{R^2}}\int_0^1 {\varepsilon
^{h-1}\hat {u}(x\varepsilon ,y\varepsilon )d\varepsilon -u_2 (x,y)} }
\right|\le $$
$$\le V_0^1  {\left( {\varepsilon ^h\hat
{u}(x\varepsilon ,y\varepsilon )} \right)}\left( {1-R^{2h}}
\right),$$

Therefore, we come to the next result.

\textbf{Theorem 4.} For the solution of the problem (\ref{eq7})-(\ref{eq9}) $u_1 ,u_2 $
approximate formulas are valid
\[
u_1 \approx \frac{1}{\ln \frac{1}{R^2}}\hat {u}_3 (x,y)-\frac{R^{2h}}{\ln
\frac{1}{R^2}}\hat {u}_3 \left( {\frac{R^2}{x},\frac{R^2}{y}} \right),
\quad
u_2 \approx \frac{1-R^{2h}}{\ln \frac{1}{R^2}}\hat {u}_3 .
\]
This theorem allows the physical interpretation: component of the solution
$u_2 $ with the accuracy to a numerical multiplier it is approximately equal
to the solution of the third homogeneous boundary problem with boundary
condition
\[
2h\hat {u}_3 +\frac{\partial }{\partial n}\hat {u}_3 =\hat {u},(x,y)\in S.
\]
It is necessary to consider the formula from [3], The formula links
solutions of the first and third homogeneous boundary problems:
\[
\hat {u}_3 =\int\limits_0^1 {\varepsilon ^{h-1}\hat {u}} (x\varepsilon
,y\varepsilon )d\varepsilon .
\]
Study the case of small values of thickness of the external layer
$1-R\approx 0$, thus we consider coefficient k$>>$1. If the expression
component $u_2 $ rewritten in the form:
\[
u_2 (x,y)=\frac{2k}{k+1}\sum\limits_{j=0}^\infty {\left( {\frac{k-1}{1+k}}
\right)} ^{2j}\hat {u}(xR^{4j},yR^{4j})-
\]
\[
-\frac{2k}{k+1}\frac{k-1}{1+k}\sum\limits_{j=0}^\infty {\left(
{\frac{k-1}{1+k}} \right)} ^{2j}\hat {u}(xR^2R^{4j},yR^2R^{4j}),
\]
define the number h using the formula
\[
\frac{k-1}{1+k}=R^{2h}.
\]
and lemma 1. We get the approximate formulas for component of the problem
solution (\ref{eq7})-(\ref{eq9}):
\[
u_1 \approx \frac{1-R^{2h}}{1+R^{2h}}\frac{1}{2\ln \frac{1}{R^2}}\left(
{\left( {\hat {u}_3 (x,y)-R^{2h}\hat {u}_3 \left( {R^2x,R^2y} \right)}
\right)+} \right),
\]
\[+R^{2h}\left( {\hat {u}_3 \left( {\frac{R^2}{x},\frac{R^2}{y}}
\right)-R^{2h}\hat {u}_3 \left( {\frac{R^4}{x},\frac{R^4}{y}} \right)}
\right)
\]
\[
u_2 \approx \frac{1-R^{2h}}{2\ln \frac{1}{R^2}}\left( {\hat {u}_3
(x,y)-R^{2h}\hat {u}_3 \left( {R^2x,R^2y} \right)} \right).
\]
\subsection{ Thin ring}

Let {\_} the function $\hat {u}=\hat {u}(x,y)$ - is harmonic in the unit
circle, $B=\{(x,y):x^2+y^2<1\}$ continuous in its closure $\bar {B}$. In the
ring $K_R $- with the smaller radius R and larger radius 1, consider the
Dirichlet problem (\ref{eq10})-(\ref{eq11}).

We apply Lemma 1, we obtain the following assessment of solutions of the
Dirichlet problem in the annulus (\ref{eq10})-(\ref{eq11})
\[
\left| {\int\limits_0^1 {\frac{\hat {u}(x\varepsilon ,y\varepsilon )-\hat
{u}\left( {\frac{R^2}{x}\varepsilon ,\frac{R^2}{y}\varepsilon }
\right)}{\varepsilon \ln \frac{1}{R^2}}} d\varepsilon -u(x,y)} \right|\le
\mathop V\limits_0^\infty ,
\]
here
\[
\mathop V\limits_0^1 =\mathop V\limits_0^1 \left( {\hat {u}(x\varepsilon
,y\varepsilon )-\hat {u}\left( {\frac{R^2}{x}\varepsilon
,\frac{R^2}{y}\varepsilon } \right)} \right)-
\]
variation of function [10] taken from a variable $\varepsilon$ on the line segment [0;
1]. Thus, the following theorem is valid.

\textbf{Theorem 5.} For the solutions of the Dirichlet problem in the
annulus (\ref{eq10})-(\ref{eq11}) approximate formula is valid
\[
u(x,y)\approx \frac{\hat {u}_2 (x,y)-\hat {u}_2 \left(
{\frac{R^2}{x},\frac{R^2}{y}} \right)}{\ln \frac{1}{R^2}},
\]
Where $\hat {u}_2 $ - the solution of the Neumann problem for the equation
Laplace's in a circle with the boundary condition
\[
\frac{\partial \hat {u}}{\partial n}=\hat {u},(x,y)\in S.
\]
\section{Conclusion}

From the results given in article follows, that the accuracy of the found
formula is order equal to the thickness of the outer layer. Therefore,
obtaining the formulas of higher accuracy class represents an important
interest. This problem is solved by using asymptotic formulas presented in
corollary 1, 2. The problem of distribution of the results on the
multilayered plates and also on the boundary problems with the more general
boundary conditions is actual.

\end{document}